\documentclass[12pt]{article}

\usepackage{amscd,amsmath, amssymb, fancyhdr} 

\numberwithin{equation}{section}

\newcommand{\version}{version 3.0,\ \  Febr 23, 2011}

\renewcommand{\leftmark}{{\scriptsize A formally K\"ahler structure on the knot space}}

\def\eqref#1{(\ref{#1})}
\newcommand{\goth}{\mathfrak}
\newcommand{\g}{{\mathfrak g}}
\newcommand{\arrow}{{\:\longrightarrow\:}}

\def\C{{\Bbb C}}
\newcommand{\R}{{\Bbb R}}

\newcommand{\6}{\partial}
\def\1{\sqrt{-1}\:}
\newcommand{\restrict}[1]{{\left|_{{\phantom{|}\!\!}_{#1}}\right.}}
\newcommand{\cntrct} {\rfloor}               

\renewcommand{\tilde}{\widetilde}
\renewcommand{\bar}{\overline}
\renewcommand{\phi}{\varphi}
\renewcommand{\epsilon}{\varepsilon}

\newcommand{\End}{\operatorname{End}}
\newcommand{\Tot}{\operatorname{Tot}}
\newcommand{\Id}{\operatorname{Id}}

\newcommand{\Alt}{\operatorname{Alt}}
\newcommand{\Lie}{\operatorname{Lie}}
\newcommand{\Diff}{\operatorname{\sf Diff}}

\newcommand{\codim}{\operatorname{codim}}

\newcommand{\hor}{{\operatorname{\sf hor}}}
\newcommand{\ver}{{\operatorname{\sf ver}}}

\newcommand{\const}{\operatorname{\sf const}}

\newcommand{\Tw}{\operatorname{Tw}}

\newcommand{\Knot}{\operatorname{Knot}}
\newcommand{\LKnot}{\operatorname{\sf LKnot}}

\newcounter{Mycounter}[section]
\newcounter{lemma}[section]
\renewcommand{\thelemma}{{Lemma \thesection.\arabic{lemma}}}
\newcommand{\lemma}{%
    \setcounter{lemma}{\value{Mycounter}}
    \refstepcounter{lemma}
    \stepcounter{Mycounter}
    {\noindent \bf \thelemma:\ }}

\newcounter{claim}[section]
\renewcommand{\theclaim}{{Claim \thesection.\arabic{claim}}}
\newcommand{\claim}{%
    \setcounter{claim}{\value{Mycounter}}
    \refstepcounter{claim}
    \stepcounter{Mycounter}
    {\noindent \bf \theclaim:\ }}

\newcounter{sublemma}[section]
\newcounter{corollary}[section]
\renewcommand{\thecorollary}{{Corollary \thesection.\arabic{corollary}}}
\newcommand{\corollary}{%
    \setcounter{corollary}{\value{Mycounter}}
    \refstepcounter{corollary}
    \stepcounter{Mycounter}
    {\noindent \bf \thecorollary:\ }}

\newcounter{theorem}[section]
\renewcommand{\thetheorem}{{Theorem \thesection.\arabic{theorem}}}
\newcommand{\theorem}{%
    \setcounter{theorem}{\value{Mycounter}}
    \refstepcounter{theorem}
    \stepcounter{Mycounter}
    {\noindent \bf \thetheorem:\ }}

\newcounter{conjecture}[section]
\newcounter{proposition}[section]
\renewcommand{\theproposition}
      {{Proposition \thesection.\arabic{proposition}}}
\newcommand{\proposition}{%
    \setcounter{proposition}{\value{Mycounter}}
    \refstepcounter{proposition}
    \stepcounter{Mycounter}
    {\noindent \bf \theproposition:\ }}

\newcounter{definition}[section]
\renewcommand{\thedefinition}
      {{Definition~\thesection.\arabic{definition}}}
\newcommand{\definition}{%
    \setcounter{definition}{\value{Mycounter}}
    \refstepcounter{definition}
    \stepcounter{Mycounter}
    {\noindent \bf \thedefinition:\ }}

\newcounter{example}[section]
\newcounter{remark}[section]
\renewcommand{\theremark}{{Remark \thesection.\arabic{remark}}}
\newcommand{\remark}{%
    \setcounter{remark}{\value{Mycounter}}
    \refstepcounter{remark}
    \stepcounter{Mycounter}
    {\noindent \bf \theremark:\ }}

\newcounter{problem}[section]
\newcounter{question}[section]
\makeatletter

\@addtoreset{equation}{section} \@addtoreset{footnote}{section}
\makeatother

\def\blacksquare{\hbox{\vrule width 5pt height 5pt depth 0pt}}
\def\endproof{\blacksquare}

\begin{document}
\begin{center}
{\LARGE\bf
A formally K\"ahler structure\\[2mm] 
on a knot space of a $G_2$-manifold\\[4mm]
}

 Misha
Verbitsky\footnote{The work is partially supported by the
 RFBR grant 10-01-93113-NCNIL-a,
RFBR grant 09-01-00242-a, Science Foundation of 
the SU-HSE award No. 10-09-0015 and AG Laboratory HSE, 
RF government grant, ag. 11.G34.31.0023\\



}

\end{center}

{\small \hspace{0.15\linewidth}
\begin{minipage}[t]{0.7\linewidth}
{\bf Abstract} \\
A knot space in a manifold $M$ is a space of 
oriented immersions $S^1 \hookrightarrow M$ 
up to $\Diff(S^1)$. J.-L. Brylinski has shown
that a knot space of a Riemannian threefold
is formally K\"ahler. We prove that a space
of knots in a holonomy $G_2$ manifold is formally
K\"ahler. 
\end{minipage}
}

\tableofcontents


\section{Introduction}


Let $M$ be an oriented Riemannian 3-fold,
and $\Knot(M)$ be its knot space, that is, a space
of non-parametrized, immersed, oriented loops,
represented by a map which is injective
outside a finite set.
J.-L. Brylinski has proved that $\Knot(M)$ is an 
infinite-dimensional formally K\"ahler Fr\'echet  
manifold (see Section
\ref{_frechet_Section_} for an explanation 
of these terms).  This formal K\"ahler structure
is easy to construct, though the proof of its
formal integrability is non-trivial. 
Given a knot $S\subset M$, its tangent space
$T_S \Knot(M)$ is a space of sections of its
normal bundle $NS$, which is 2-dimensional,
oriented and orthogonal. A 2-dimensional
oriented Euclidean vector space has a natural
complex structure, which is defined through
counter-clockwise turns. Therefore,
the bundle $NS$ is a 1-dimensional
complex Hermitian bundle. Therefore,
the space of sections of $NS$ is a complex Hermitian
Fr\'echet vector space. The corresponding 
Hermitian form is easy to obtain from the
volume 3-form on $M$ by integration along
the knots (\ref{_fiberwi_integra_Definition_}).

$G_2$-manifolds appear naturally
as a main object of ``octonionic
algebraic geometry'', playing
the same role for octonions as 
hyperk\"ahler and hypercomplex manifolds
play for quaternions. The main engine
for the study of quaternionic geometry
is the twistor construction, which makes
a complex manifold from a manifold with
a quaternionic structure. It is well known
that the twistor data can be used 
to reconstruct the quaternionic structure.
Singularities in hyperk\"ahler and hypercomplex
geometry and many natural geometric objects
can also be studied in terms of twistors
(\cite{_Verbitsky:Symplectic_II_}, 
\cite{_Verbitsky:hypercomple_}).

One would expect the hypothetical
octonion twistor space (if it exists) to
bring similar benefits. However, none
of the usual approaches to constructing
complex structures on twistor manifolds 
works for $G_2$-geometry, and it seems
that something must be sacrificed.
In the present paper, we sacrifice 
finite-dimensionality of a twistor space.

We propose a twistor-like
construction resulting in a formally K\"ahler
structure on the knot space of a $G_2$-manifold
(\ref{_knot_complex_main_Theorem_}). 
This construction is similar in flavour
to one of J.-L. Brylinski; in fact, our
approach to the proof of formal integrability
is similar to the argument
of L. Lempert (\cite{_Lempert:CR_}),
who used a CR twistor space constructed
for $G_2$-manifolds by C. LeBrun.
A $G_2$-analogue of LeBrun's twistor
space was constructed in \cite{_Verbitsky:CR-twistor_},
and now we use it to study the complex
structure on the knot space. We also
interpret several objects of $G_2$-geometry
(instanton bundles, associative subvarieties)
as holomorphic objects on the knot space
(Section \ref{_G_2_geome_Section_}).

The symplectic structure which appears in 
this construction was previously obtained
by M. Movshev (\cite{_Movshev_}).


\section{Fr\'echet manifolds and formally K\"ahler geometry}
\label{_frechet_Section_}


In this section, we briefly introduce Fr\'echet manifolds
and basic geometric structures on such manifolds. For
a detailed exposition, please see \cite{_Lempert:Dolbeault_}.

\subsection{Fr\'echet manifolds and knot spaces}

Recall that a {\bf Fr\'echet space} is an infinite-dimensional
topological vector space $V$ admitting a translation-invariant
complete metric. It is equivalent to say that $V$ has a
countable family of seminorms $\| \cdot\|_1$, 
$\|\cdot\|_2$, $\| \cdot\|_3$, ...,  and the topology
on $V$ is induced by a complete, translation-invariant 
metric
\[
d(x,y):= \sum_{i=1}^\infty \frac 1 {2^i} \min(\|x-y\|_i,1).
\]
(\cite{_Bourbaki:EVT5_}).
{\bf A differentiable map} of Fr\'echet spaces is
a map which can be approximated at each point by a continuous 
linear map, up to a term which decays faster than linear,
in the sense of this metric. In a similar way one defines
smooth (infinitely differentiable) maps of Fr\'echet spaces.

A {\bf Fr\'echet manifold} is a ringed space, locally
modeled on a space of differentiable functions on a
Fr\'echet space, with transition functions smooth.

When $M$ is a compact finite-dimensional manifold,
the space $C^\infty M$ of smooth functions on $M$
has a sequence of norms $C_k$, with
\[
\| f\|_{C^k}= \sum_{i=0}^k\sup_M |f^{(i)}|,
\]
where $f^{(i)}$ denotes the $i$-th derivative.
It is well known that this system of seminorms
is complete on $C^\infty(M)$, giving a 
structure of Fr\'echet space on $C^\infty(M)$.
Similar constructions allow one to define the 
Fr\'echet structure on the space of smooth 
sections of a vector bundle.

This is used to define a structure of a Fr\'echet
manifold on various infinite-dimensional spaces
arising in geometry, in particular on a space
$\operatorname{Imm}(X,M)$ of smooth immersions $X\hookrightarrow M$, 
and on a group of diffeomorphisms, which becomes
a Fr\'echet Lie group. 

The quotient $\operatorname{Imm}(X,M)/\Diff(X)$
is a Fr\'echet orbifold, locally modeled on the 
total space $NX$ of its normal bundle. To see this,
one needs to construct a slice of the $\Diff(X)$-action,
which can be done easily using a Riemannian metric.
The orbifold points correspond to those maps
which are wrapped several times on themselves.
Denote by $\operatorname{Imm}_0(X,M)$
the space of immersions which are injective
outside of a positive codimension set of
self-intersections.  Clearly, then
$\operatorname{Imm}_0(X,M)/\Diff(X)$ is 
a Fr\'echet  manifold.

For the present paper, the most important Fr\'echet
manifold is a space $\Knot(M):= \operatorname{Imm}_0(S^1,M)/\Diff_+(S^1)$
of oriented knots (non-parametrized immersed loops,
injective outside of a finite set) in $M$.
We could work with the orbifold
$\operatorname{Imm}(S^1,M)/\Diff_+(S^1)$
instead, and all the results would remain
valid in the orbifold context. 
To simplify terminology, we work with manifolds and restrict 
ourselves to $\operatorname{Imm}_0(S^1,M)/\Diff_+(S^1)$.

\subsection{Formally complex Fr\'echet manifolds}
\label{_formally_co_Subsection_}

Let $F$ be a Fr\'echet manifold. One can define 
{\bf the sheaf of vector fields} $TF$ on $F$ as a shief
of continuous derivations of its structure sheaf. 
A commutator of two derivations is again a derivation. 
This gives a Lie algebra structure on the sheaf
of vector fields. 

The formally integrable almost complex structures
are defined as usual, in the following way.

\hfill

\definition
Let $F$ be a Fr\'echet manifold, and $I:\; TF \arrow TF$
a smooth $C^\infty F$-linear endomorphism of the tangent
bundle satisfying $I^2=-1$. Then $I$ is called
{\bf an almost complex structure on $F$}.

\hfill

\remark
Clearly, $I$ defines a decomposition
$TF\otimes \C= T^{1,0}F\oplus T^{0,1} F$,
where $T^{1,0}F$ is the $\1$-eigenspace of $I$, and 
$T^{0,1}F$ the $-\1$-eigenspace. Indeed, 
$x = \frac 1 2 (x + \1 Ix) + \frac 1 2 (x - \1 Ix)$.

\hfill

\definition
An almost complex structure on a Fr\'echet manifold $(F,I)$ is called
{\bf formally integrable},  if 
$[T^{1,0}F, T^{1,0}F] \subset T^{1,0}F$,
where $[\cdot,\cdot]$ denotes the commutator of vector
fields. In this case $(F,I)$ is called {\bf a formally
  complex manifold}. 

\hfill

\remark
Just as it happens in the finite-dimensional case, the projection
of $[T^{1,0}F, T^{1,0}F]$ to $T^{0,1}F$ is always
$C^\infty F$-linear. This gives an operation
\[ \Lambda^2 T^{1,0}F \stackrel N \arrow T^{0,1}F,\] called
{\bf the Nijenhuis tensor}. The Nijenhuis tensor
of an almost complex Fr\'echet manifold $(F,I)$ 
vanishes if and only if it is  formally integrable.

\hfill

\definition
A function $f$ on an almost complex Fr\'echet 
manifold is called {\bf holomorphic} if
$\langle df, X\rangle=0$ for any vector field
$X \in T^{0,1}F$. A smooth map to a complex topological vector
space is called {\bf holomorphic} if its composition
with continuous complex linear functionals is always holomorphic.

\hfill

\definition
An almost complex structure on a Fr\'echet manifold $F$ is called
{\bf strongly integrable} if there exists an atlas of local
coordinate charts which are given by holomorphic maps
to complex Fr\'echet spaces. In this case $F$ is called
{\bf holomorphic}.

\hfill

\remark
Every holomorphic Fr\'echet manifold is 
formally integrable, which is obvious, because
one can locally generate $T^{1,0}F$ by coordinate vector
fields which commute. The converse implication is known
to be false. For finite-dimensional manifolds, formal integrability
implies integrability of a complex structure, by
a deep analytic result called the Newlander-Nirenberg 
theorem (\cite{_Newla_Nire:integra_}). An infinite-dimensional version of 
Newlander-Nirenberg theorem is false (see \cite{_Lempert:CR_}).

\hfill

\definition
Let $(F,I)$ be a formally integrable almost complex
Fr\'echet manifold, $g$ a Hermitian structure on $F$, and
$\omega$ be the corresponding $(1,1)$-form. We say that
$(F,I,g)$ is {\bf formally K\"ahler} if $\omega$ is closed.

\subsection{Formally K\"ahler structures on knot spaces}

Let $M$ be a smooth manifold, and $\Knot(M)$ be its knot space.
As we have mentioned already, $\Knot(M)$ is a Fr\'echet
manifold. Locally at $l\in \Knot(M)$, this manifold
is modeled on the space of smooth sections of a normal
bundle $Nl$. 

Geometric structures on the space of knots on an oriented
3-manifold $M^3$ were a subject of much research
(see e.g. \cite{_Brylinski:preprint_},
\cite{_Le_Brun:sheets_}, \cite{_Lempert:CR_},
and the book \cite{_Brylinski_}). In
\cite{_Brylinski:preprint_},
a formally K\"ahler structure on $\Knot(M^3)$
was constructed.  In \cite{_Lempert:CR_}, it was shown that this
formally complex structure is never strongly integrable.

In his book \cite{_Brylinski_}, J.-L. Brylinski
gives many uses for the formal K\"ahler structure
on the space of knots. Another possible application 
of the formally K\"ahler structure on the space of knots
(not explored much, so far) is to the spaces of discriminants
of knots and their cohomology. V. A. Vassiliev
defined the eponymous knot invariants by considering
the stratification on the space of knots by 
successive discriminant spaces. Later, M. Kontsevich
redefined some of these cohomology spaces and proved
that they carry the mixed Hodge structure.
It is easy to see that the discriminant spaces
are in fact complex subvarieties, in the sense
of formally complex structure on the knots. 
One would expect that the mixed Hodge structure
on Vassiliev invariants comes from this
complex stratification, just as it would
happen in the finite-dimensional case. 

The aim of this paper is to generalize these results
to $G_2$-manifolds, which are 7-dimensional Riemannian
manifolds with special holonomy group which lies 
in $G_2$.

The formally complex structure on the
space of knots of a 3-dimensional manifold $M$
can be defined in terms of the vector product 
on $TM$. Indeed, let $S$ be a knot in $M$, and
$\gamma'$ a unit tangent vector field to $S$.
A vector product with $\gamma'$ defines a
complex structure on the normal bundle $NS$,
which is used to define the formal
complex structure on $\Knot(M)$.
In \cite{_Lee_Leung:knots_} and \cite{_Lee_Leung:insta_},
geometry of manifolds with vector products was
explored at some length, and many results
about the knot spaces and instantons were
obtained from a similar vector product
construction.


\section{$G_2$-manifolds}


\subsection{$G_2$-geometry: basic notions}

$G_2$-manifolds originally appeared
in Berger's classification of Riemannian holonomy
(\cite{_Berger:list_}, \cite{_Besse:Einst_Manifo_}).
The first examples of $G_2$-manifolds were
obtained by R. Bryant and S. Salamon (\cite{_Bryant_Salamon_}).
The compact examples of $G_2$-manifolds
were constructed by D. Joyce (\cite{_Joyce_G2_}, 
\cite{_Joyce_Book_}). 
In this introduction we follow the approach
to $G_2$-geometry which is due to
N. Hitchin (see \cite{_Hitchin:3-forms_}).

\hfill

\definition
Let $\rho\in \Lambda^2 \R^7$ be a 3-form on $\R^7$.
We say that $\rho$ is {\bf non-degenerate} if the 
dimension of its stabilizer is maximal:
\[
\dim St_{GL(7)}\rho = \dim GL(7) - \dim \Lambda^3(\R^7) =
49-35 =14.
\]
In this case, $St(\rho)$ is one of two real forms 
of a 14-dimensional Lie group $G_2(\C)$. 
We say that $\rho$ is {\bf non-split} if it satisfies
$St(\rho|_x)\cong G_2$, where $G_2$ denotes the
compact real form of $G_2(\C)$.
{\bf A $G_2$-structure} on a 7-manifold
is a 3-form $\rho \in \Lambda^3(M)$, which is
non-degenerate and non-split at each point $x\in M$. We shall always
consider a $G_2$-manifold as a Riemannian
manifold, with the Riemannian structure induced
by the $G_2$-structure as follows. 

\hfill

\remark \label{_metric_volume_G_2_Remark_}
A form $\rho$ defines
a $\Lambda^7 M$-valued metric on $M$:
\begin{equation}\label{_metric_vol_valued_Equation_}
g(x,y) = (\rho\cntrct x)\wedge (\rho \cntrct y) \wedge \rho
\end{equation}
(we denote by $\rho \cntrct x$ the contraction of $\rho$
with a vector field $x$). The Riemannian volume form
associated with this metric gives a section of 
$\Lambda^7 M\otimes (\Lambda^7 M)^{7/2}$. Squaring and
taking the 9-th degree root, we obtain a trivialization of the 
volume. Then \eqref{_metric_vol_valued_Equation_} defines
a metric $g$ on $M$, by construction $G_2$-invariant.

\hfill

\definition 
An  $G_2$-structure is called {\bf an integrable $G_2$-\-struc\-ture},
if $\rho$ is preserved by the corresponding Levi-Civita connection.
An integrable 
$G_2$-manifold is a manifold equipped with an integrable $G_2$-structure.
Holonomy group of such a manifold clearly lies in $G_2$;
for this reason, the integrable $G_2$-manifolds are often called
{\bf holonomy $G_2$-manifolds}. 

\hfill

\remark
In the literature, ``the $G_2$-manifold'' often means
a ``holonomy $G_2$-manifold'', and ``$G_2$-structure''
``an integrable $G_2$-structure''. A $G_2$-structure
which is not necessarily integrable is called 
``an almost $G_2$-structure'', taking analogy
from almost complex structures. Further on in this
paper, we shall follow this usage, unless 
specified otherwise.

\hfill

\remark
As shown in \cite{_Fernandez_Gray:G_2_},
integrability of a $G_2$-structure induced
by a 3-form $\rho$ is equivalent to $d\rho = d(*\rho)=0$.
For this reason the 4-form $*\rho$ is called
{\bf the fundamental 4-form of a $G_2$-manifold},
and $\rho$ {\bf the fundamental 3-form}.

\hfill

\remark\label{_SU(3)_Remark_}
Let $V=\R^7$ be a 7-dimensional real space equipped
with a non-degenerate 3-form $\rho$ with $St_{GL(7)}(\rho)=G_2$.
As in \ref{_metric_volume_G_2_Remark_}, one can easily
see that $V$ has a natural $G_2$-invariant metric. 
For each vector $x\in V$, $|x|=1$, its stabilizer 
$St_{G_2}(x)$ in $G_2$ is isomorphic to $SU(3)$.
Indeed, the orthogonal complement $x^\bot$ is
equipped with a symplectic form $\rho \cntrct x$,
which gives a complex structure $g^{-1} \circ (\rho \cntrct x)$
as usual. This gives an embedding $St_{G_2}(x)\hookrightarrow U(3)$.
Since the space of such $x$ is $S^6$, and
the action of $G_2$ in $S^6$ is transitive, one has
$\dim St_{G_2}(x)= \dim G_2 - \dim S^6 =8 = \dim U(3)-1$.
To see that $St_{G_2}(x)= SU(3)\subset U(3)$ and not some
other codimension 1 subgroup, one should notice that
$St_{G_2}(x)$ preserves two 3-forms $\rho \restrict{x^\bot}$
and $\rho^*\cntrct x \restrict{x^\bot}$, where $\rho^* = *\rho$
is the fundamental 4-form of $V$. A simple linear-algebraic
calculation implies that $\rho \restrict{x^\bot}+\1 
\rho^*\cntrct x \restrict{x^\bot}$ is a holomorphic
volume form on $x^\bot$, which is clearly preserved by 
$St_{G_2}(x)$. Therefore, the natural embedding
$St_{G_2}(x)\hookrightarrow U(3)$ lands $St_{G_2}(x)$
to $SU(3)$. Using the dimension count $\dim St_{G_2}(x)=\dim SU(3)$
(see above), we show that the embedding
$St_{G_2}(x)\hookrightarrow SU(3)$ is also surjective.

\subsection{Octonion structure and a vector product}
\label{_Octo_Subsection_}

Let $V=\R^7$ be a 7-dimensional space equipped
with a non-degenerate, non-split constant 3-form $\rho$ 
inducing a $G_2$-action on $V$. Then $V$ is equipped with
the {\bf vector product}, defined as follows: 
$x\star y = \rho(x, y, \cdot)^\sharp$. Here
$\rho(x, y, \cdot)$ is a 1-form obtained by
contraction, and $\rho(x, y, \cdot)^\sharp$
is its dual vector field.

\hfill

\remark\label{_comple_str_vector_pr_Remark_}
The complex structure
on an orthogonal complement $v^\bot$ is 
given by a vector product: $x\arrow v\star x$,
if $|v|=1$.

\hfill

It is not hard to see
that $(V, \star)$ becomes isomorphic to the
imaginary part of the octonion algebra,
with $\star$ corresponding to half of the commutant.
In fact, this is one of the many ways used to define
an octonion algebra. The whole octonion algebra
is obtained as ${\Bbb O}:=V \oplus \R$,
with the product given by 
\[ (x, t) (y, t')= 
   (ty + t'x + x\star y, g(x,y)+tt')
\]
Here, $x, y$ and $ty + t'x + x\star y$
are vectors in $V$, and $t, t', g(x,y)+tt'\in \R$.

Given two non-collinear vectors in $V$, they generate
a quaternion subalgebra in octonions. When these two vectors
satisfy $|v|=|v'|=1$, $v\bot v'$, the standard basis $I, J, K$
in imaginary quaternions can be given by a triple 
$v, v', v\star v'\in V$.

A 3-dimensional subspace $A \subset V$ is called
{\bf associative} if it is closed under the vector
product. The set of associative subspaces is in
bijective correspondence with the set of quaternionic
subalgebras in octonions.

\subsection{A CR twistor space of a $G_2$-manifold}

\definition
Let $M$ be a smooth manifold, $B\subset TM$ be a subbundle
in its tangent bundle, and $I\in \End B$ be an automorphism,
satisfying $I^2=-\Id_B$. Consider the (1,0) and (0,1)-bundles
$B^{1,0}, B^{0,1}\subset B \otimes \C$, which are the
eigenspaces of $I$ corresponding to the eigenvalues
$\1$ and $-\1$. The sub-bundle $B^{1,0}\subset TM \otimes \C$ 
is called {\bf a CR-structure on $M$} if it is {\em involutive},
that is, it satisfies $[B^{1,0}, B^{1,0}]\subset B^{1,0}$.

\hfill

Let $M$ be an almost $G_2$-manifold. 
From \ref{_SU(3)_Remark_} it follows that with every
vector $x\in TM$, $|x|=1$, one can associate a
complex Hermitian structure on its orthogonal
complement $x^\bot$. The easiest way to define 
this structure is to notice that $x^\bot$
is equipped with a symplectic structure 
$\rho\cntrct x$ and a metric $g\restrict{x^\bot}$,
which can be considered as real and imaginary
parts of a complex-valued Hermitian 
product. Then the complex structure is obtained as
usual, as $I:=(\rho\cntrct x)\circ g^{-1}$.

\hfill

\corollary\label{_rho_and_Herm_Corollary_}
Let $(M,\rho)$ be an almost $G_2$ manifold, $m\in M$ a
point, and $x\in T_m M$ a non-zero vector. Then
 the symplectic form $\rho\cntrct x$
is a Hermitian form of a natural complex
structure on $x^\bot \subset T_m M$.
\endproof

\hfill

\definition\label{_B^1,0_Definition_}
Consider the unit sphere bundle $S^6M$
over $M$, with the fiber $S^6$, and let $T_{hor} S^6M \subset TS^6 M$
be the horizontal sub-bundle corresponding to the
Levi-Civita connection. This sub-bundle has a natural
section $\theta$; at each point $(x, m)\in S^6 M$,
$m\in M$, $x\in T_m M$, $|x|=1$, we take
$\theta\restrict{(x, m)}= x$, using the standard 
identification $T_\hor S^6 M\restrict{(x, m)}= T_m M$.
Denote by $B\subset T_\hor S^6 M$ the orthogonal complement to $\theta$
in $T_{hor} S^6 M$. Since at each point $(x, m)\in S^6 M$,
the restriction $B \restrict{(x, m)}$ is identified
with $x^\bot \subset T_m M$, this bundle is equipped
with a natural complex structure, that is, an operator
$I \in \End B$, $I^2= -\Id_B$. 

\hfill

\theorem\label{_CR_twi_Theorem_}
(\cite{_Verbitsky:CR-twistor_})
Let $M$ be an almost $G_2$-manifold, $S^6M\subset TM$
its unit sphere bundle, and $B\subset T S^6 M$
a sub-bundle of its tangent bundle constucted
above, and equipped with the complex structure $I$.
Then $B^{0,1}\subset B\otimes \C \subset T S^6 M\otimes \C$
is involutive if and only if $M$ is a holonomy $G_2$-manifold. 
\endproof

\hfill

\definition\label{_twistor_Definition_}
Let $M$ be a holonomy $G_2$-manifold,
and \[ \Tw(M):=(S^6 M, B, I)\] the CR-manifold constructed 
in \ref{_CR_twi_Theorem_}. Then $\Tw(M)$ is called
{\bf a CR-twistor space of $M$}.

\hfill

The key argument in the proof of \ref{_CR_twi_Theorem_}
is achieved by constructing a certain 3-form
$\Omega\in \Lambda^3(S^6 M, \C)$. This form is
of type $(3,0)$ on $B$ with respect to the complex 
structure constructed in \ref{_B^1,0_Definition_}, and satisfies 
$d\Omega\restrict B=0$. The restriction $\Omega\restrict B$
is determined uniquely from the $SU(3)$-structure on $B$:
we construct $\Omega$ in such a way that $\Omega\restrict B$
is equal to the holomorphic volume form associated with
the $SU(3)$-structure. Since we are going to use
this form and the expression for $d\Omega$ obtained in 
\cite{_Verbitsky:CR-twistor_}, we describe it explicitly below.

Consider the fundamental 3-form and 4-form $\rho$ and $\rho^*:= *(\rho)$
on a holonomy $G_2$-manifold 
$M$, and let $\theta\in T_\hor S^6 M$ be the tautological
vector field constructed in \ref{_B^1,0_Definition_}.
Denote by $\pi:\; S^6 M \arrow M$ the standard projection.
Then 
\begin{equation}\label{_Omega_defi_Equation_}
\Omega:= \pi^* \rho + \1 (\pi^*\rho^*\cntrct \theta),
\end{equation}
where $\pi^*\rho^*\cntrct \theta$ denotes the contraction
of $\pi^*\rho^*$ and $\theta \in T_\hor S^6 M\subset TS^6M$.

\hfill

The 4-form $d\Omega:=\1d (\pi^*\rho^*\cntrct \theta)$
was computed in \cite{_Verbitsky:CR-twistor_} explicitly,
as follows. 

Consider a natural embedding 
\begin{equation}\label{_S^6_to_Lambda^3_Equation_}
S^6 M \stackrel\phi\hookrightarrow \Tot(\Lambda^3 M)
\end{equation}
into the total space of $\Lambda^3 M$, mapping 
$(v,m)$ to $(\pi^*\rho^*\cntrct v, m)$.
Let $\Xi$ be a 4-form on 
$\Tot(\Lambda^3 M)$ written in local coordinates 
$p_1, ..., p_7$ as 
\begin{equation}\label{_Xi_defi_Equation_}
\Xi = 
\sum_{i_1 < i_2< i_3} dq_{i_1, i_2, i_3}
\wedge  dp_{i_1}\wedge dp_{i_2}\wedge dp_{i_3}
\end{equation}
where $q_{i_1, i_2, i_3}$ is a function
on $\Tot(\Lambda^3 M)$, linear on fibers and expressed
as 
\[ 
   q_{i_1, i_2, i_3}:= 
   \frac d {dp_{i_1}}\wedge \frac d {dp_{i_3}} \wedge \frac d {dp_{i_3}}
\]
(here we identify the 3-vector fields on $M$ with 
linear functions on $\Tot(\Lambda^3 M)$).
The form $\Xi$ is a 4-dimensional analogue
of the usual Hamiltonian 2-form on $\Tot(\Lambda^1 M)$,
and satisfies similar standard properties. It is
called ``the Hamiltonian 4-form'' in \cite{_Verbitsky:CR-twistor_}.

In  \cite{_Verbitsky:CR-twistor_} the following 
result was proven.

\hfill

\proposition\label{_d_Omega_Proposition_}
Let $M$ be a holonomy $G_2$-manifold, and
$\Omega$ be the 3-form on $S^6 M$ constructed above.
Then $-\1 d\Omega= d (\pi^*\rho^*\cntrct \theta)$
is equal to $\phi^* \Xi$, where 
$S^6 M \stackrel\phi\hookrightarrow \Tot(\Lambda^3 M)$
is the embedding defined in \eqref{_S^6_to_Lambda^3_Equation_}.
\endproof


\section{A K\"ahler structure on a knot space
of a $G_2$-manifold}


\subsection{Knot spaces on CR-manifolds}

\def\SS{{\Bbb S}}
\newcommand{\Imm}{\operatorname{Imm}}

Let $\SS$ and $M$ be smooth manifolds, 
$\Imm(\SS, M)$ the set of immersions
from $\SS$ to $M$, and $\Knot(\SS, M):= \Imm(\SS,M)/\Diff(\SS)$
the corresponding knot orbifold. Clearly,
$\Knot(\SS, M)$ is a Fr\'echet orbifold, modeled
on the space of sections of a normal bundle
$NS$ in a neighbourhood of a point $S\in \Knot(\SS, M)$.

\hfill

\definition
Suppose that $(M,B,I)$ is a CR-manifold. A knot 
$S\in \Knot(\SS, M)$ is called {\bf transversal},
if for all $s\in S$, the intersection
of $T_s S\cap B\restrict s=0$, that is,
$B$ is transversal to $TS$ everywhere.
Denote the space of transversal knots by
$\Knot_B(\SS, M)$.

\hfill

\remark
Let $(M,B,I)$ be a CR-manifold,
and $\Knot_B(\SS,M)$ be the space of transversal knots.
Suppose that $\dim \SS = \codim B$.
For each knot $S\in \Knot_B(\SS,M)$, one has
$NS\cong B\restrict S$, hence $\Knot_B(\SS,M)$
has a natural almost complex structure.

\hfill

\theorem\label{_transv_knots_integra_Theorem_}
Let $(M,B,I)$ be a CR-manifold, $\SS$ be a smooth manifold 
for which $\dim \SS = \codim B$,
and $\Knot_B^0(\SS,M)$ the space of transversal knots
which are embedded to $M$ (that is, have no self-\-in\-ter\-sec\-tion),
equipped with a complex structure as above.
Then $\Knot_B^0(\SS,M)$ is a formally
complex Fr\'echet manifold.

\hfill

{\bf Proof:} 
Consider an embedded knot $S\in \Knot_B^0(\SS,M)$. The space
$T^{1,0}_S \Knot_B^0(\SS,M)$ of $(1,0)$-tangent vectors
is by definition equal to the space of sections 
of $B^{1,0}\restrict S$. Let $X,Y\in B^{1,0}\restrict S$
be some sections of the bundle $B^{1,0}\restrict S$.
We extend $X,Y$ to sections $X_1, Y_1$ of $B^{1,0}$
in some neighbourhood of $S$ (this is possible, because
$S$ has no self-intersections). The vector fields
$X_1, Y_1\in TM$ are used in a usual way to define
vector fields $\tilde X, \tilde Y\in T^{1,0} \Knot_B^0(\SS,M)$
satisfying $\tilde X\restrict S= X$, $\tilde Y\restrict S=Y$.
Denote by $Z_1=[X_1, Y_1]$ the commutator of $X_1, Y_1$.
Since $(M,B,I)$ is a CR-manifold, $Z_1\in B^{1,0}$.
Denote by $\tilde Z\in T^{1,0} \Knot_B^0(\SS,M)$ the corresponding
vector field on $\Knot_B^0(\SS,M)$, defined in a neighbourhood
of $S$. Clearly, $[\tilde X, \tilde Y]= \tilde Z$, and
its $(0,1)$-component vanishes.
Therefore, the Nijenhuis tensor 
\[ 
  N:\; T^{1,0}_S \Knot_B^0(\SS,M)\times T^{1,0}_S \Knot_B^0(\SS,M) \arrow
 T^{0,1}_S \Knot_B^0(\SS,M)
\]
vanishes on arbitrarily chosen vectors $X,Y$. We proved that
$N=0$, hence $\Knot_B^0(\SS,M)$ is integrable. \endproof

\hfill

The same argument also proves the following theorem.

\hfill

\theorem\label{_integra_cr_on_knots_Theorem_}
Let $(M,B,I)$ be a CR-manifold, and $F\subset TM$ 
be a sub-bundle of $TM$ containing $B$. Consider the space
$\Knot^0_F(\SS, M)$, defined as above, and let
${\cal B} \subset T \Knot^0_F(\SS, M)$ be a sub-bundle
consisting of all $\xi\in T_S \Knot^0_F(\SS, M)$
which belong to 
$\Gamma(B\restrict S) \subset \Gamma(F\restrict S)=T_S \Knot^0_F(\SS, M)$
(here $\Gamma$ denotes the space of global sections of a vector bundle).
Consider a complex structure operator on $\cal B$ induced
by $I:\; B\restrict S \arrow B\restrict S$.
Then $(B,I)$ is an integrable CR-structure on $\Knot^0_F(\SS, M)$.

\endproof

\subsection{Two failed proofs of formal integrability and a definition 
of $L$-knots}

It looks natural to obtain the integrability of the almost
complex structure on $\Knot(M)$ for a $G_2$-manifold $M$ directly
from \ref{_integra_cr_on_knots_Theorem_}. The first naive
approach is to use the standard projection 
$\Knot^0_F(S^6M)\stackrel\pi\arrow \Knot(M)$,
where $\Knot^0_F(S^6M)=\Knot^0_F(S^1, S^6M)$ is the space of 
knots transversal to the bundle $F=T_\ver S^6M\oplus B$.
It is easy to see that the projection of a knot in $S^6M$
transversal to $T_\ver S^6M$ is a knot in $M$.
By \ref{_integra_cr_on_knots_Theorem_} and 
\ref{_CR_twi_Theorem_}, $\Knot^0_F(S^6M)$ is equipped
with an integrable CR-structure ${\cal B}^{1,0}$
induced from $B^{1,0}\subset TS^6M \otimes \C$. The projection
$\Knot^0_F(S^6M)\stackrel\pi\arrow \Knot(M)$ induces a
map
\begin{equation}\label{_proje_on_tan_Equation_}
d\pi:\; {\cal B} \arrow T\Knot(M)
\end{equation}
which is an isomorphism at each point of $S^6M$
(it follows directly from an isomorphism $T_S \Knot(M)=\Gamma(B\restrict S)$,
where $\Gamma(B\restrict S)$ is the space of sections of 
the restriction $B\restrict S$). Were $d\pi$ also 
compatible with the complex structures, integrability of
the complex structure on $\Knot(M)$ would follow 
immediately from the integrability of ${\cal B}^{1,0}$.
Unfortunately, this is not so. In fact, it is easy to
see that the projection \eqref{_proje_on_tan_Equation_}
is compatible with the complex structures only 
at the one point in each fiber of $\pi$ 
(\ref{_L_knot_only_where_complex_lineat_Remark_}). 
This one point corresponds to a special class of knots
in $S^6M$ called {\em Legendrian} in \cite{_Lempert:CR_}.
To avoid confusion with the usual Legendrian knots
on the contact manifold $S^6M$, we call the
Lempert's Legendrian knots
{\bf the $L$-knots}.

\hfill

\definition
Let $M$ be a Riemannian 7-manifold, 
$S\in \Knot_{T_\ver S^6M}(S^6M)$ a transversal knot on $S^6M$,
and $S_1\in \Knot(M)$ its projection to $M$.
Consider a unit speed parametrization $\gamma_1:\; S^1 \arrow S_1$,
and let $\gamma:\; S^1 \arrow S$ be the corresponding
parametrization in $S$. Each point $t\in S^1$ gives
a unit vector $\dot\gamma_1(t)\in T_{\gamma_1(t)}M$.
Assume that $\gamma(t)=(\dot\gamma_1(t), \gamma_1(t))$.
Then $S$ is called {\bf an $L$-knot}.

\hfill

We return now to the case when $M$ is a $G_2$-manifold. 

\hfill

\remark\label{_L_Knot_section_Remark_}
Clearly, there is precisely one $L$-knot in every
fiber of the projection 
$\Knot^0_{T_\ver S^6M}(S^6M)\stackrel\pi\arrow \Knot(M)$.
This gives a smooth section $\Psi:\; \Knot(M)\arrow \LKnot(M)$,
where $\LKnot(M)\subset \Knot^0_F(S^6M)$ denotes the
space of $L$-knots.

\hfill

\remark\label{_L_knot_only_where_complex_lineat_Remark_}
For each $L$-knot $S$, the map 
\begin{equation}\label{_d_proje_at_Lknot_Equation_}
d\pi:\; {\cal B}\restrict S \arrow T_{\pi(S)}\Knot(M)
\end{equation}
is complex linear. Indeed, for each $(v,m)\in S$,
the complex structure in $B\restrict S$
is given by a vector product with a vector
$v$, and the complex structure in the normal
bundle $N\pi(S)$
is given by a vector product with  a unit
tangent vector to $\pi(S)$. These vectors
are equal precisely when $S$ is an $L$-knot. 

From this argument, it is clear that
\eqref{_d_proje_at_Lknot_Equation_} is complex
linear only if $S$ is an $L$-knot.

\hfill

This leads us to the second failed argument
for the proof of integrability of $\Knot(M)$.
It is natural to expect the embedding 
\[ \Psi:\; \Knot(M)\arrow \LKnot(M)\subset \Knot_F(S^6M) \]
together with \ref{_CR_twi_Theorem_}
to bring us integrability of the
complex structure on $\Knot(M)$.
Unfortunately, this argument does
not work, because the bundle ${\cal B}\subset T\Knot_F(S^6M)$
is not tangent to the submanifold $\LKnot(M)\subset \Knot_F(S^6M)$.

This is why a direct application of 
\ref{_CR_twi_Theorem_} to the integrability of complex
structure on $\Knot(M)$ does not work. Instead, we use
a different approach, which employs \ref{_CR_twi_Theorem_}
only as a framework. 

\subsection{The tangent space to the space of $L$-knots}

To prove the integrability of the knot space, we use an
explicit description of the tangent space $T\LKnot(M)$.
The following claim is trivial.

\hfill

\claim\label{_diffeo_on_LKnot_Claim_}
Consider the natural action of the group $\Diff(M)$
on the space $S^6M$, interpreted as a double cover of
the projectivization ${\Bbb P}TM$. Then 
$\LKnot(M)\subset \Knot(S^6M)$ is a  $\Diff(M)$-invariant
subset of the corresponding knot space.

\hfill

{\bf Proof:} The definition of $\LKnot(M)$ is functorial 
with respect to diffeomorphisms, hence $\LKnot(M)$
is clearly $\Diff(M)$-invariant. \endproof

\hfill

Since the natural projection $\pi:\; \LKnot(M)\arrow \Knot(M)$
is a diffeomorphism of Fr\'echet manifolds, to describe the
tangent space $T\LKnot(M)\subset T\Knot(S^6M)$, we need to
describe the image of vector fields $X\in T\Knot(M)$
in $T\LKnot(M)\subset T\Knot(S^6M)\restrict{\LKnot(M)}$.

\hfill

\proposition\label{_tangent_to_LKnot_Proposition_}
We consider knot spaces on a Riemannian manifold $M$
of dimension $7$.\footnote{In fact, the proof of 
\ref{_tangent_to_LKnot_Proposition_} would
hold for any dimension $>2$.}
Let $S_1\in \Knot(M)$, and let $S\in \LKnot(M)$ be the
corresponding $L$-knot. Consider the decomposition
\begin{equation}\label{_tangent_to_Lkno_decompo_Equation_}
  T_S\Knot(S^6M)= \Gamma(T_\ver (S^6M)\restrict S) 
  \oplus \Gamma(B\restrict S)
\end{equation}
obtained from the Levi-Civita-induced orthogonal 
decomposition
$TS^6M= T_\ver (S^6M)\oplus B\oplus \R \cdot\theta$
(\ref{_B^1,0_Definition_}).
Let $\gamma_1$ be a unit speed parametrization of $S_1$.
Consider a vector $X\in T_S\LKnot(M)\subset T_S\Knot(S^6M)$, and
let $X_1\in T_{S_1}\Knot(M)=\Gamma(B\restrict S)$ be the corresponding
tangent vector to $S_1$. Then the decomposition
of $X$ into components of \eqref{_tangent_to_Lkno_decompo_Equation_}
is written as
\begin{equation}\label{_decompo_tange_via_conne_Equation_}
X = -\nabla_{\dot\gamma_1} X_1 + X_1,
\end{equation}
where $\nabla_{\dot\gamma_1} X_1$ is a 
$\Gamma(T_\ver (S^6M)\restrict S)$-component of $X$, and 
$X_1$ the $\Gamma(B\restrict S)$-component.

\hfill

{\bf Proof:} Clearly, it would suffice to prove
\eqref{_decompo_tange_via_conne_Equation_} for a dense open
subset $\Knot^0(M)\subset \Knot(M)$ of embedded knots. For such a knot,
we may extend $X_1$ to a vector field $\tilde X_1\in TM$.
Let $\tilde X$ be the correspoding vector field
tangent to $\LKnot(M)$ (\ref{_diffeo_on_LKnot_Claim_}).
Clearly, $\tilde X\restrict S = X$. The action of
$\tilde X$ on $\LKnot(M)$ is a restriction of a 
vector field $\tilde X_{S^6M}\in T\Knot(S^6M)$ acting on
$S^6M$ by functoriality (here we again
identify $S^6M$ with a double cover of ${\Bbb P}TM$).
To finish the proof of \ref{_tangent_to_LKnot_Proposition_}
it would remain to compute the vector field $\tilde X_{S^6M}$.

At a point $(\dot\gamma_1, m)\in S\subset S^6M$, 
the vector field $\tilde X_{S^6M}$ is equal to 
$(-\nabla_{\dot\gamma_1} X_1, X_1)$, which gives
\begin{equation}\label{_action_on_S^6M_Equation_}
\tilde X_{S^6M}\restrict{S}=(-\nabla_{\dot\gamma_1} X_1, X_1).
\end{equation}
Indeed, for any section $V\in TM$, 
\[
\Lie_{\tilde X} V = [\tilde X, V] = \nabla_{\tilde X} V - \nabla_V\tilde X,
\]
and choosing $V$ in such a way that $\nabla_{\tilde X} V=0$, and $V=v$,
we obtain that
\[
e^{t\tilde X}(v,m)=(v-t\nabla_v \tilde X+ o(t^2), e^{t\tilde X}m).
\]
Then \eqref{_action_on_S^6M_Equation_} clearly follows.
From \eqref{_action_on_S^6M_Equation_}, 
\eqref{_decompo_tange_via_conne_Equation_} is apparent.
\endproof

\hfill

For the proof of integrability of the complex structure on 
$\Knot(M)$, the following proposition is used.

\hfill

\proposition\label{_Xi_on_LKnot_Proposition_}
Let $M$ be a holonomy $G_2$-manifold, and 
$\Xi\in \Lambda^4(S^6M)$ the 4-form obtained as in 
\eqref{_Xi_defi_Equation_}. Consider the corresponding
4-form $\tilde \Xi$ on $\Knot(S^6M)$ mapping 
vector fields $X_1, ..., X_4\in T_S \Knot(S^6M)$
to the integral $\int_\gamma \Xi(X_1, ..., X_4)dt$,
where $\gamma(t)$ is a unit speed parametrization of $S$.
Then $\tilde \Xi\restrict{\LKnot(M)}=0$.

\hfill

{\bf Proof:} Let $X_i^\ver$, $X_i^\hor$
be the vertical and horizontal components
of $X_1\in T_S\Knot(S^6M)$, under the natural
decomposition 
\[
  T_S\Knot(S^6M)= \Gamma(T_\ver (S^6M)\restrict S) 
  \oplus \Gamma(B\restrict S).
\]
In \cite{_Verbitsky:CR-twistor_}, it was shown that
$\Xi\restrict{T_\hor S^6M}$ vanishes. Also, from the
local formula \eqref{_Xi_defi_Equation_} it is apparent
that $\Xi(v_1, v_2, \cdot, \cdot)=0$, where
$v_1, v_2 \in T_\ver (S^6M)$. Therefore,
\[
\Xi(X_1, ..., X_4) = \Alt\rho^*(X_1^\ver, X_2^\hor, X_3^\hor, X_4^\hor),
\]
where $\Alt$ denotes the skew-symmetrization over the indexes 
$1,2,3,4$. However, by \ref{_tangent_to_LKnot_Proposition_},
$X_i^\ver = -\nabla_{\dot \gamma}X_i^\hor$, hence
\begin{equation}\label{_tilde_Xi_expli_Equation_}
\tilde\Xi(X_1, ..., X_4) = \int_{\gamma} 
\Alt\rho^*(\nabla_{\dot\gamma}X_1^\hor, X_2^\hor, X_3^\hor, X_4^\hor).
\end{equation}
Since $M$ is a holonomy $G_2$-manifold, 
 $\nabla_{\dot\gamma}\rho^*=0$, and \eqref{_tilde_Xi_expli_Equation_}
gives
\[
\tilde\Xi(X_1, ..., X_4) = 
\int_\gamma \frac d{dt}\left(\rho^*(X_1^\hor\restrict{\gamma(t)}, ..., X_4^\hor\restrict{\gamma(t)}\right)dt,
\]
which is equal 0 by Stokes' theorem. \endproof

\subsection{Non-degenerate $(3,0)$-forms}

\definition
Let $M$ be a smooth or Fr\'echet manifold, equipped 
with an almost complex structure, and $\Omega\in \Lambda^{3,0}(M)$ 
a $(3,0)$-form. We say that $\Omega$ is {\bf non-degenerate}
if for any $X\in T^{1,0}(M)$ there exist $Y,Z\in T^{1,0}(M)$
such that $\Omega(X,Y,Z)\neq 0$.

\hfill

The utility of non-degenerate (3,0)-forms is due to the following
simple theorem.

\hfill

\theorem\label{_3-form_closed_the_integra_Theorem_}
Let $M$ be a smooth or a Fr\'echet manifold equipped 
with an almost complex structure, and $\Omega\in \Lambda^{3,0}(M)$ 
a non-degenerate $(3,0)$-form. Assume that $d\Omega=0$.
Then $M$ is formally integrable.

\hfill

{\bf Proof:} Let $X,Y \in T^{1,0}(M)$ and $Z,T\in T^{0,1}(M)$.
Since $\Omega$ is a (3,0)-form, it vanishes on 
$(0,1)$-vectors. Then Cartan's formula together with 
$d\Omega=0$ implies that
\begin{equation}\label{_Cartan's_Eqution_}
0=d\Omega(X,Y,Z,T)= \Omega(X,Y,[Z,T]).
\end{equation}
From the non-degeneracy of $\Omega$ we obtain that
unless $[Z,T]\in T^{0,1}(M)$, for some $X,Y \in T^{1,0}M$
one would have $\Omega(X,Y,[Z,T])\neq 0$. 
Therefore, \eqref{_Cartan's_Eqution_} implies
$[Z,T]\in T^{0,1}(M)$, for all $Z,T\in T^{0,1}(M)$,
which means that $M$ is integrable. 
\endproof

\hfill

For the proof of integrability of the almost complex
structure on the knot space, we also need the following
trivial lemma.

\hfill

\lemma\label{_tilde_Omega_nondege_on_knots_Lemma_}
Let $M$ be a $G_2$-manifold, and 
$\Omega:= \pi^*\rho +\1(\pi^*\rho^*\cntrct\theta)\in \Lambda^3(S^6M)$
be the 3-form constructed in \ref{_Omega_defi_Equation_}.
Denote by $\tilde \Omega$ the corresponding form on
$\LKnot(M)$, mapping $X_1,X_2,X_3\in T_S \LKnot(S^6M)\subset T_S \Knot(S^6M)$
to $\int_{\gamma(t)} \Omega(X_1,X_2,X_3) dt$, where
$\gamma(t)$ is a unit speed parametrization of $S$. 
Then $\Omega$ a non-degenerate $(3,0)$-form on $\LKnot(M)\cong \Knot(M)$.

\hfill

{\bf Proof:} The space $T_S \Knot(M)$ is identified with
the space of sections $\Gamma(B\restrict S)$, where
$B$ is an $SU(3)$-bundle defined in \ref{_B^1,0_Definition_}.
From its construction, it is clear that 
$\Omega\restrict {T_\ver S^6M}$ vanishes, and
on $B$ the form $\Omega$ is equal to the
standard complex volume form.
Therefore, $\Omega$ is of type $(3,0)$,
and for each $X\in B^{1,0}\restrict S$, there
exist $Y,Z\in B^{1,0}\restrict S$ such that
$\Omega(X,Y,Z)$ is everywhere real, non-negative,
and positive at any point where $X\neq 0$.
Then, $\tilde\Omega(X,Y,Z)$ is real and positive,
unless $X=0$. \endproof

\subsection{The proof of integrability of the knot space}
\label{_integra_Subsection_}

Let $M$ be a holonomy $G_2$-manifold, and
$S\subset M$ {\em a knot}, that is, a class of an immersion
$S^1\stackrel \gamma \hookrightarrow M$ up to 
oriented reparametrizations, and injective outside of a
finite set. We can assume
that $S$ is parametrized, with $|\gamma'| = \const$
(such a parametrization is obviously unique).

The tangent space $T_S \Knot(M)$ is identified
with the space of sections of the normal bundle $NS$.
At each point $s\in S$, $N_sS= (T_sS)^\bot$
is the orthogonal complement to an oriented line
$T_sS \subset T_s M$. Then \ref{_SU(3)_Remark_}
gives a complex structure on $N_s S$. This defines
an almost complex structure on $\Knot(M)$. To define
a Hermitian form, the following construction is used.

\hfill

\definition\label{_fiberwi_integra_Definition_}
Let $\Knot^m(M)\subset \Knot(M) \times M$ be the space
of {\em marked knots}, that is, pairs 
$(S^1\stackrel\gamma\hookrightarrow M, s\in S^1)$,
where $\gamma$ is injective somewhere, and $|\gamma'|=1$. Clearly, 
the forgetful map $\Knot^m(M)\stackrel\pi \arrow \Knot(M)$
is an $S^1$-fibration. The fiberwise integration map
\[
\Lambda^i(\Knot^m(M))\stackrel {\pi_*}\arrow 
\Lambda^{i-1}(\Knot(M))
\]
is defined as usual, 
\[
   \pi_*(\alpha)\restrict S := \int_{\pi^{-1}(S)}
   \left(\alpha\cntrct \frac d {dt}\right) dt
\]
where $t$ is a parameter on $S$.
It is easy to check that $\pi_*$ commutes with the de Rham 
differential. Define $\sigma:\; \Knot ^m (M) \arrow M$
as follows, 
\[ 
  \sigma(S^1\stackrel\gamma\hookrightarrow \Knot(M), s\in
  S^1):= \gamma(s).
\]
This gives an interesting map
\[
\pi_*\sigma^*:\; \Lambda^i(M) \arrow 
\Lambda^{i-1}(\Knot(M))
\]
commuting with the de Rham differential. 

For a $G_2$-manifold $(M, \rho)$, the 2-form
$\pi_* \sigma^*(\rho)$ was computed by M. Movshev
in \cite{_Movshev_}, who proved that it is symplectic.

\hfill

\claim\label{_push_of_rho_Claim_}
Let $(M, \rho)$ be an almost $G_2$-manifold, $S\in
\Knot(M)$ a knot, and $\alpha, \beta \in NS$ two sections
of a normal bundle, considered as tangent vectors
$a,b \in T_S\Knot(M)$. Consider the 
integral $S(a,b):=\int_S \rho(a,b, \cdot)\restrict S$.
Then $\pi_* \sigma^*(\rho)(a,b) = S(a,b)$.

\hfill

{\bf Proof:} This claim is essentially a restatement of a definition
(see \cite{_Movshev_} for more detail).
\endproof

\hfill

Comparing \ref{_push_of_rho_Claim_}
and \ref{_rho_and_Herm_Corollary_}, 
we obtain the following result.

\hfill

\proposition
Let $(M, \rho)$ be an almost $G_2$-manifold
$\omega:= \pi_* \sigma^*(\rho)$ the Movshev's
2-form on $\Knot(M)$, and $I$ the almost complex structure
on $\Knot(M)$ constructed above.  Then
$(\Knot(M), I, \omega)$ is an almost
complex Hermitian Fr\'echet manifold.

\endproof

\hfill

The main result of this paper is the 
following theorem.

\hfill

\theorem\label{_knot_complex_main_Theorem_}
Let $M$ be a holonomy $G_2$-manifold, and 
$(\Knot(M), I, \omega)$ an almost
complex Hermitian Fr\'echet manifold
constructed above. Then $(\Knot(M), I, \omega)$
is formally K\"ahler. 

\hfill

\remark
The manifold $(\Knot(M),  \omega)$ is symplectic
(\cite{_Movshev_}). This is clear from the
construction of $\omega=\pi_* \sigma^*(\rho)$,
because $\pi_* \sigma^*$ commutes with the de
Rham differential.

\hfill

{\bf Proof of \ref{_knot_complex_main_Theorem_}:}
To prove integrability of the almost complex structure
on $\Knot(M)$, we identify $\Knot(M)$ with the space
of $L$-knots $\LKnot(M)\subset \Knot(S^6M)$ as in
\ref{_L_Knot_section_Remark_}. As follows
from \ref{_tilde_Omega_nondege_on_knots_Lemma_},  the space $\LKnot(M)$
is equipped with a non-degenerate $(3,0)$-form
$\tilde\Omega$. From \ref{_d_Omega_Proposition_} it follows
that $d\tilde\Omega=\1\tilde\Xi$, and from
\ref{_Xi_on_LKnot_Proposition_} that $\tilde \Xi\restrict{\LKnot(M)}=0$,
hence the form $\tilde\Omega$ is closed on $\LKnot(M)$.
Now, integrability of the almost complex structure
on $\Knot(M)$ follows from \ref{_3-form_closed_the_integra_Theorem_}.
\endproof


\section{The complex structure on the knot space and $G_2$-geometry}
\label{_G_2_geome_Section_}


\subsection{Associative subvarieties of a $G_2$-manifold
and complex subvarieties in  its knot space}

The complex geometry of a knot space can be used to study
the geometry of a $G_2$-manifold. Many notions of a 
$G_2$-geometry can be directly translated to the
language of complex geometry, as follows.

\hfill

\definition
Let $X\subset M$ be a 3-dimensional subvariety
of a $G_2$-manifold. We say that $X$ is {\bf associative}
if $T_xX\subset T_xM$ is an associative subspace
for each smooth point $x\in X$ (see 
Subsection \ref{_Octo_Subsection_} for a definition
of an associative subspace).

\hfill

\proposition\label{_curves_in_knot_Proposition_}
Let $M$ be a holonomy $G_2$-manifold, and $A\subset \Knot(M)$
a 1-dimensional complex subvariety. Denote 
by $\tilde A\subset M$ the union of all knots in $A$.
Then $\tilde A$ is an associative subvariety of $M$.

\hfill

{\bf Proof:} 
Let $\gamma\in A$ be a knot, and $x, y \in T_\gamma A$ two 
tangent vectors, considered as  sections of a normal bundle $N\gamma$,
with $I(x)=y$. A complex structure on $T_\gamma A$ is given
by a vector product with the unit vector field 
$\frac{\gamma'}{|\gamma'|}$ 
(\ref{_comple_str_vector_pr_Remark_}). Therefore, 
the 3-dimensional space $\langle x, y, \gamma'\rangle$
is closed under the vector product. \endproof

\hfill

\proposition
Let $M$ be a holonomy $G_2$-manifold,
and $X \subset M$ a subvariety, $1< \dim X <7$.
Then $\Knot(X) \subset \Knot(M)$ 
is a formally complex subvariety if and only if 
$X$ is an associative subvariety.

\hfill

{\bf Proof:} The same argument as 
in \ref{_curves_in_knot_Proposition_}
proves that $T_xX \subset T_xM$ is closed with respect
to the vector product, for any smooth point $x\in X$.
However, any proper subalgebra of octonions
is isomorphic to $\R$, $\C$ or ${\Bbb H}$,
as an easy algebraic argument implies.
Since  $1< \dim X <7$, this is a quaternion
subalgebra, and the subspace $T_x X$ 
3-dimensional and associative. 
\endproof

\hfill

\subsection{Holomorphic bundles on a knot space}

$G_2$ instanton bundles were introduced in 
\cite{_Donaldson_Thomas_}, and much studied since then.
This notion is a special case of a more general
notion of an instanton on a calibrated manifold,
which is already well developed. Many estimates 
known for 4-dimensional manifolds
(such as Uhlenbeck's compactness theorem)
can be generalized to the calibrated case 
(\cite{_Tian:Calibrated_}, \cite{_Tian_Tao_}).

Recently, $G_2$-instantons became a focus of
much activity because of attempts to construct
a higher-dimensional topological quantum field theory,
associated with $G_2$ and 3-dimensional Calabi-Yau
manifolds (\cite{_Donaldson_Segal_}).

\hfill

\definition
Let $M$ be a $G_2$-manifold, and 
$\Lambda^2 M = \Lambda^2_7(M) \oplus \Lambda^2_{14}(M)$
the irreducible decomposition of the bundle of 2-forms 
$\Lambda^2(M)$ associated with the $G_2$-action.
A vector bundle $(B, \nabla)$ with connection
is called {\bf a $G_2$-instanton} if its curvature
lies in $\Lambda^2_{14}(M)\otimes \End(B)$.

\hfill

\remark
Since the curvature of a holonomy $G_2$-manifold
lies in $\Lambda^2M \otimes \g_2$, and $g_2\subset \goth{so}(TM)$
is identified with $\Lambda^2_{14}$ under the identification 
$\goth{so}(TM)=\Lambda^2(M)$, the curvature of $TM$ lies
in $\Lambda^2_{14}M \otimes\End(TM)$. Therefore, 
the tangent bundle and all its tensor powers are 
$G_2$-instantons. 

\hfill

\remark
Let $M$ be a finite-dimensional complex manifold,
and $E$ a Hermitian bundle on $M$. Recall that a holomorphic
structure on $E$ induces a unique Hermitian connection $\nabla$
on $E$ which its curvature $\Theta$ satisfying 
$\Theta \in \Lambda^{1,1}(M)\otimes \End E$ and 
$\nabla^{0,1}=\bar\6$, where $\bar\6$ is the holomorphic
structure operator
(this connection is called {\bf the Chern connection}).
This motivates the following definition.

\hfill

\definition
Let $(F,I)$ be a formally complex Fr\'echet manifold,
and $(E,\nabla)$ a Hermitian bundle with connection.
We say that $(E, \nabla)$ is {\bf formally holomorphic}
if the curvature $\Theta$ of $\nabla$ satisfies
$\Theta \in \Lambda^{1,1}(F)\otimes \End E$.

\hfill

\remark
Let $E$ be a vector bundle with connection
on a Riemannian manifold $M$, and $\Knot(M)$
be the knot space of $M$. For a given $S\in \Knot(M)$,
consider the space $E(S)$ of sections of $E\restrict S$.
Consider an infinite-dimensional bundle $\tilde E$ on $\Knot(M)$
with fiber $E(S)$ at $S\in \Knot(M)$. This bundle
can be obtained as $\pi_\bullet\sigma^* E$,
where $\sigma:\; \Knot ^m (M) \arrow M$,
$\pi:\; \Knot ^m (M) \arrow \Knot(M)$
are the maps defined in Subsection 
\ref{_integra_Subsection_}, and $\pi_\bullet$
is a pushforward, considered in the sense
of sheaf theory. Also, every
connection on $E$ induces a connection
$\tilde \nabla:= \pi_\bullet \sigma^* \nabla$ on $\tilde E$.

\hfill

\theorem
Let $M$ be a $G_2$-manifold,
$\Knot(M)$ its knot space equi\-p\-ped with
a natural formally K\"ahler structure, $(E, \nabla)$
a Hermitian vector bundle with connection,
and $(\tilde E, \tilde \nabla)$ the corresponding
bundle on $\Knot(M)$. Then $(\tilde E, \tilde \nabla)$ 
is formally holomorphic if and only if $\nabla$
is a $G_2$-instanton.

\hfill

{\bf Proof:}
Clearly, the curvature $\tilde \Theta$ of $\tilde E$ is obtained
by lifting the curvature $\Theta$ of $E$ to $\Knot(M)$
in a natural way. From \cite[Proposition 3.2]{_Verbitsky:CR-twistor_},
it follows that a form belongs to $\Lambda^2_{14}(M)$
if and only if its restriction to each 6-dimensional
subspace $x^\bot \subset T_m M$ is of type $(1,1)$.
This is equivalent to $\tilde \Theta$ being
of type (1,1) on $\Knot(M)$. \endproof

\hfill

\hfill

\noindent{\bf Acknowledgements:} 
I am grateful to L\'aszl\'o Lempert for his interest
and for reminding me about the orbifold points
in the space of non-parametrized immersions, to Stefan 
Nemirovsky for finding an error in the proof, and
to Michael Entov for his consultations on Legendrian knots.

{
\small

\noindent {\sc Misha Verbitsky\\
{\sc Laboratory of Algebraic Geometry, SU-HSE,\\
7 Vavilova Str. Moscow, Russia, 117312}\\
\tt verbit@maths.gla.ac.uk, \ \  verbit@mccme.ru\\
}}

\end{document}